%%%%%This a plain TeX file of the paper
%%%%%Three Recitations ... Series
%%by Doron Zeilberger
 
%%possibly adjust :
\hoffset=1cm
\voffset=1cm
\magnification=1200
\hsize=11.25cm
\vsize=18cm
\parskip 0pt
\parindent=12pt
\frenchspacing

\def\qed{\quad\raise -2pt\hbox{\vrule\vbox to 10pt{\hrule 
width 4pt \vfill\hrule}\vrule}}

\def\halmos{\penalty 500 \hbox{\qed}\par\smallskip}

\font\eightbf=cmbx8
\font\eightrm=cmr8 \font\sixrm=cmr6 
\font\eighttt=cmtt8

\long\def\proclaim #1. #2\endproclaim{\medbreak
{\bf #1.\enspace}{\sl#2}\par\medbreak}

\catcode`\@=11

\def\eightpoint{%
  \textfont0=\eightrm \scriptfont0=\sixrm 
\scriptscriptfont0=\fiverm
  \def\rm{\fam\z@\eightrm}%
   \abovedisplayskip=9pt plus 2pt minus 6pt
  \abovedisplayshortskip=0pt plus 2pt
  \belowdisplayskip=9pt plus 2pt minus 6pt
  \belowdisplayshortskip=5pt plus 2pt minus 3pt
  \smallskipamount=2pt plus 1pt minus 1pt
  \medskipamount=4pt plus 2pt minus 1pt
  \bigskipamount=9pt plus 3pt minus 3pt
  \normalbaselineskip=9pt
  \setbox\strutbox=\hbox{\vrule height7pt depth2pt width0pt}%
  \let\bigf@ntpc=\eightrm \let\smallf@ntpc=\sixrm
  \normalbaselines\rm}

\def\appeln@te{}
\def\vfootnote#1{\def\@parameter{#1}\insert
  \footins\bgroup\eightpoint
  \interlinepenalty\interfootnotelinepenalty
  \splittopskip\ht\strutbox %
  \splitmaxdepth\dp\strutbox \floatingpenalty\@MM
  \leftskip\z@skip \rightskip\z@skip
  \ifx\appeln@te\@parameter\indent \else{\noindent #1\ }\fi
  \footstrut\futurelet\next\fo@t}

\def\footnoterule{\kern-6\p@
  \hrule width 2truein \kern 5.6\p@} % the \hrule is .4pt high

\catcode`\@=12

%begin macros

\def\CF{\hbox{\eightrm CF}}
\def\gcd{\mathop{\rm g.c.d.}}
\def\deg{\mathop{\rm deg}}
\font\eightrm=cmr8 \font\sixrm=cmr6 
\font\eighttt=cmtt8
%end macros

{\eightpoint
\noindent
S\'eminaire Lotharingien de Combinatoire, B24a, 1990

}
\vskip 1.5cm
\centerline{\bf THREE RECITATIONS ON HOLONOMIC SYSTEMS}
\medskip
\centerline{\bf AND HYPERGEOMETRIC SERIES}
\bigskip
\centerline{D{\sevenrm ORON} ZEILBERGER \footnote{$^1$}
{\eightpoint
Department of Mathematics, Temple University,
Philadelphia, PA 19122, USA.\hfil\break
{\eighttt (zeilberg@math.temple.edu).} 
Supported in part by the NSF.\hfil\break
Updated version, February 22, 1995, of a paper that appeared in
Publ. I.R.M.A. Strasbourg, 1993, 461/S-34, Actes $24^e$
S\'eminaire Lothringien, p. 5-37.}}
 
\bigskip\medskip
{\eightpoint
{\eightbf Preface:} These ``recitations'' were given in the
$24^{th}$ session of the ``S\'eminaire Lotharingien'', held in
the Spring of 1990, somewhere in the Vosges mountains.  I thank
Dominique Foata for inviting me, and letting me sample one of
these charming seminars that preserve the spirit that
Oberwolfach lost a long time ago. I would like to thank Peter
Paule and Volker Strehl for the invitation to include them in
this special issue of the JSC, and for many helpful comments.
}

\bigskip
\centerline{\bf Foreword} 

\medskip
When we teach calculus we have lectures and
recitations. These notes are meant as ``recitations'' or
something like ``Schaum outlines'' for the theory. The role of the
``lectures'' or ``textbook'' is provided by Gosper's
path-breaking paper ``{\it A Decision Procedure for Indefinite
summation}'', Proc. Nat. Acad. Sci. USA {\bf 75} (1978), 40-42, and
by the following papers by myself and my  collaborators, Gert
Almkvist and Herb Wilf.

\smallskip
[AZ] (With Gert Almkvist) {\it The method of differentiating under the
integral sign}, J. Symbolic Computation {\bf 10}, 571-591 (1990).
 
[WZ1] (With H. S. Wilf) {\it Rational functions certify combinatorial
identities}, J. Amer. Math. Soc. {\bf 3}, 147-158 (1990).
 
[WZ2] (With H. S. Wilf) {\it Towards computerized proofs of identities},
Bulletin of the Amer. Math. Soc. {\bf 23}, 77-83 (1990).

[Z1] {\it A Holonomic systems approach to special functions
identities}, J. of Computational and Applied Math. {\bf 32},
321-368 (1990).
 
[Z2] {\it  A Fast Algorithm for proving terminating hypergeometric
identities}, Discrete Math {\bf 80}, 207-211 (1990).
 
[Z3]{\it The method of creative telescoping}, J. Symbolic Computation
{\bf 11}, 195-204 (1991).

\medskip
In addition, the following papers give further expositions by
myself.

\smallskip 
[Z5] {\it Identities in search of identity}, 
J. Theoretical Computer Science {\bf 117},  23-38 (1993).  
 
[Z6] {\it Theorems for a price: Tomorrow's semi-rigorous
mathematical culture}, Notices of the Amer. Math. Soc. {\bf 40 \#
8} (Oct. 1993), 978-981.\hfil\break
Reprinted (followed by a
critique by George Andrews) in: Math. Intell. {\bf 16 \#1} 11-14.

\medskip
There also appeared superb expositions by Pierre Cartier~[C], on
the general theory, and by Tom Koornwinder~[K], on the fast
algorithm and its q-analog. Excellent treatments of Gosper's
algorithm  and of the fast algorithm
are given in  sections 5.7 and 5.8 of~[GPK] below.
The former section also appears in the first edition, the latter
section is new to the second edition. More recently, Herb
Wilf~[W] wrote beautiful lecture notes.

\smallskip
[C] P. Cartier, {\it D\'emonstration ``automatique'' d'identit\'es
et fonctions hyperg\'eom\'etriques [d'apres D. Zeilberger]},
S\'eminaire Bourbaki, expos\'e $n^{o}$ $746$, Ast\'erisque {\bf
206}, 41-91, SMF, 1992.
 
[GKP] R. L. Graham, D. E. Knuth, and O. Patashnik, {\it Concrete Mathematics},
Second Edition, Addison-Wesley, Reading, 1993.
 
[K] T. H. Koornwinder, {\it Zeilberger's algorithm and its q-analogue},
J. of Computational and Applied Math. {\bf 48}, 91-111 (1993).
 
[W] H.S. Wilf, {\it ``Identities and their computer proofs"}, 
{\eightrm SPICE}
lecture notes {\bf 31}, 1993. Available by anonymous ftp
to  {\tt ftp.cis.upenn.edu} as file {\tt pub/wilf/lecnotes.ps}.

\medskip
The following papers offer important extensions, implementations,
and applications.

\smallskip 
[PS] P. Paule and M. Schorn, {\it A Mathematica version of Zeilberger's
algorithm for proving binomial coefficient identities}, J. Symbolic
Comp., to appear.
 
[Ko1] W. Koepf, {\it {\eightrm REDUCE} package for the indefinite
and definite summation }, Konrad-Zuse-Zentrum f\"ur
Informationstechnik Berlin ({\eightrm ZIB}),\hfil\break  Technical Report
TR 94-9, 1994.
 
[Ko2] W. Koepf, {\it Algorithms for the Indefinite and Definite 
Summation}, Konrad-Zuse-Zentrum f\"ur Informationstechnik
Berlin  ({\eightrm ZIB}), Technical Report TR 94-33,
1994.
 
[St1] V. Strehl, {\it Binomial sums and identities}, Maple Technical
Newsletter {\bf 10}, 37-49 (1993).

\bigskip

\centerline{\bf Recitation I: Elimination}
 
\medskip
The process of {\it elimination} consists of getting simple, or
desirable, equations out of a given system of equations.
For example
$$
(i)\ 2x \, + \, 3y \, -5 \, =0, \qquad
(ii)\ 3x-y-2=0.
$$
In order to eliminate $x$, we do
$$
3(i)-2(ii)=6x+9y-15-6x+2y+4=0
$$
getting $11y-11=0$ and hence $y=1$.

The {\it resultant} of two polynomials $P(x)$ and $Q(x)$ is
obtained by eliminating $x$ between them. The vanishing of the
resultant is the condition that they have
a common root. For example, if $f= a x^2 + b x + c$ and
$g= a' x^2 + b' x +c'$, then we have 
$$
\pmatrix{
a & b & c  &  0 \cr
0 & a  & b &  c \cr
a' & b' & c' & 0 \cr
0 & a' & b' & c'  \cr
}
\pmatrix{ {x^3} \cr
      { x^2 } \cr
      {x} \cr
      {1} \cr }
 \, = \, 
\pmatrix{ {0} \cr
      {0} \cr
      {0} \cr
      {0} \cr } .
$$
Eliminating $x$, one gets the determinant of the above matrix 
(the so-called {\it Sylvester matrix}), which is the resultant. 
 
The {\it discriminant} of a polynomial $P(x)$ is the 
resultant of $P(x)$ and $P'(x)$, and its vanishing gives the 
condition that it has a double root.
For example, for the generic second degree polynomial,
$P(x)= a x^2 +  b x + c$, eliminating $x$ from
$P(x)$ and $P'(x)$ yields
$$
-4aP(x)+(b+2ax)P'(x)= b^2 - 4ac.
$$
 
For systems of polynomial equations with several variables
$$
P_1 ( x_1 , \ldots , x_n ) = 0 , \ldots ,
P_m ( x_1 , \ldots , x_n ) = 0 ,
$$
we can eliminate $m-1$ variables,
getting a polynomial equation
$$
Q( x_1 , \ldots , x_{n-m+1} ) \, = \, 0.
$$
{\eightrm BUCHBERGER'S AMAZING GR\"OBNER BASES DO THAT FAST.}
 
\bigskip
\centerline{\bf The Joy of Operator Notation}
 
\medskip
Let $N$ be the shift operator in $n$ : $N f(n) := f(n+1)$.

\medskip
\noindent
{\it  Example:} Prove that
$$
F_{n+4}= F_{n+2} + 2 F_{n+1} + F_n,
$$
where $F_n$ are the Fibonacci numbers.
 
\medskip
\noindent 
{\it  Verbose Proof:}
$$\eqalignno{
F_{n+2} - F_{n+1} - F_n &=0,
&(i)\cr
F_{n+3} - F_{n+2} - F_{n+1} &=0,
&(ii)\cr
F_{n+4} - F_{n+3} - F_{n+2}&=0;
&(iii)\cr
F_{n+4} - F_{n+2} -2 F_{n+1} - F_n 
&=0.&(i)+(ii)+(iii)\cr}
$$

\noindent
{\it  Terse Proof}:
$$\displaylines{\qquad
( N^2 - N - 1 ) F_n=0
\Rightarrow ( N^2 + N +1) ( N^2 - N
-1) F_n =0  \hfill\cr
\hfill{}  \Rightarrow  ( N^4 - N^2 -2N-1)
F_n  =  0. \quad\cr}
$$

\goodbreak
If a sequence satisfies one recurrence, then it satisfies an
infinite number of other recurrences:
$$
P(N,n) a(n) =  0    \Rightarrow   [Q(N,n)P(N,n)] a(n)  =  0
$$
for every operator $Q(N,n)$.

In two variables, $(n,k)$, we introduce the shift operators
$N$, $K$ acting on discrete functions $F(n,k)$, by
$$
N F(n,k):= F(n+1,k),\qquad K F(n,k) := F(n,k+1).
$$
For example, the Pascal triangle equality
$$\displaylines{
 { {n+1} \choose  {k+1} } =
 { {n} \choose  {k+1} } 
 \,  +  \,   { {n} \choose  {k} }\cr 
\noalign{\hbox{is written, in operator notation, as}}
(NK-K-1) { {n} \choose {k} }  \, = \, 0.\cr}
$$
If a discrete function $F(n,k)$ satisfies two partial linear
recurrences 
$$\displaylines{
P(N,K,n,k) F(n,k) =  0 ,\qquad Q(N,K,n,k)F(n,k)=0,\cr
\noalign{\hbox{then it satisfies many, many others:}}
\{ A(N,K,n,k) P(N,K,n,k) + B(N,K,n,k)
Q(N,K,n,k) \} F(n,k)  =  0,\cr}
$$
where $A$ and $B$ can be any linear partial recurrence
operators. 

So far, everything was true for arbitrary linear
recurrence operators. From now on we will only allow linear
recurrence operators  {\it with polynomial coefficients}. The set
of linear recurrence operators with polynomial coefficients,
denoted by $C<n,k,N,K>$ is a (non-commutative) associative
algebra  generated by $N,K,n,k$ subject to the relations
$NK=KN$, $Nk=kN$, $Kn=nK$, $nk=kn$, $Nn=(n+1)N$, $Kk=(k+1)K$. 
By a clever choice of the operators $A$ and $B$, we can get the
operator in the braces above, call it $R(N,K,n)$, to be independent
of $k$.

Now write
$R(N,K,n)=S(N,n)+(K - 1) \overline R  (N,K,n)$,
where $S(N,n):=R(N,1,n)$.  Since
$R(N,K,n) F(n,k) \equiv 0$, we have
$$
S(N,n) F(n,k) = (K-1) [ - \overline R (N,K,n) F(n,k) ].
$$
Calling the function inside the square brackets above $G(n,k)$,
we get
$$
S(N,n) F(n,k) = (K-1) G(n,k).
$$
Note that if $F(n, \pm \infty ) =0$ for every $n$, then the same
is true of $G(n, \pm \infty )$. Now summing the above w.r.t. $k$
yields
$$
S(N,n) \Bigl( \sum_{k} F(n,k)\Bigr) =
\sum_{k} \bigl(G(n,k+1)  - G(n,k) \bigr ) =  0 .
$$
So $a(n):= \sum\limits_{k} F(n,k)$
satisfies the recurrence
$S(N,n) a(n) = 0 $.

\goodbreak
\medskip
\noindent
{\it  Example:}
$$\eqalignno{
F(n,k)&= {n!\over (k!\,(n-k)!},\cr
{{F(n+1,k)} \over {F(n,k)}} &= {{n+1} \over {n-k+1}},&(i)\cr
{{F(n,k+1)} \over {F(n,k)}}& = {{n-k} \over {k+1}} .&(ii)\cr}
$$
Cross multiply:
$$\eqalignno{
(n-k+1)F(n+1,k)-(n+1)F(n,k)&=0,&(i)\cr
(k+1)F(n,k+1)-(n-k)F(n,k)&=0,&(ii)\cr}
$$
In operator notation,
$$\displaylines{\noalign{\vskip-8pt}
(i) \  [(n-k+1)N-(n+1)]F \equiv 0 ,\quad
(ii) \  [(k+1)K-(n-k)]F \equiv 0 .\cr
\noalign{\hbox{Expressing the operators 
in descending powers of
$k$, we get}}
(i)\   [ (-N)k \, + \, (n+1)N-(n+1) ] F \equiv 0,\quad 
(ii) \   [(K+1) k  \, - \,  n ] F \equiv 0.\cr
\noalign{\hbox{Eliminating $k$, we get}}
(K+1)(i)+N(ii)=
\{ (K+1) [(n+1)N-(n+1)]+N(-n) \} F \equiv 0 ,\cr
\noalign{\hbox{which becomes}}
(n+1)[NK-K-1] F \equiv 0  .\cr
\noalign{\hbox{So we got that}}
R(N,K,n)=(n+1)[NK-K-1],\quad
S(N,n)=R(N,1,n)=(n+1)[N-2],\cr
\noalign{\hbox{and so we have proved the deep result that}}
a(n):= \sum_{k} { {n} \choose {k}}\cr
\noalign{\hbox{satisfies}}
(n+1)(N-2)a(n) \equiv 0  ,\cr}
$$
i.e., in everyday notation, $(n+1)[a(n+1)-2a(n)] \equiv 0$, and
hence, since $a(0) = 1$, we get that $a(n)= 2^n $.

\bigskip
\noindent 
{\it  Important observation of Gert Almkvist:} So far we had two
stages:
$$\eqalignno{
R(N,K,n)&=A(N,K,n,k)P(N,K,n,k)+B(N,K,n,k)Q(N,K,n,k)
,\qquad&(i)\cr
R(N,K,n)&=S(N,n)+(K-1) \overline R (N,K,n),&(ii)\cr
S(N,n)&=AP+BQ+(K-1) (- \overline R ),\cr}
$$
where $\overline R$  has the nice but {\it superfluous} property of
not involving $k$! {\eightrm WHAT A WASTE}. So we are led to
formulate the following.
 
\goodbreak
\bigskip
\noindent
{\it  Modified Elimination Problem:}

Input: Linear partial recurrence operators with polynomial
coefficients $P(N,K,n,k)$ and $Q(N,K,n,k)$. Find operators $A,B,C$
such that
$$
S(N,n):=AP+BQ+(K-1)C
$$
does not involve $K$ and $k$. 
 
\medskip
\noindent
{\it  Remark.} Note something strange: we are allowed
to multiply $P$ and $Q$ by any operator {\it from the left} ,
but not from the right, while we are allowed to multiply $K-1$
by any operator {\it from the right}, but not from the left.
In other words we have to find a non-zero operator, depending
on $n$ and $N$ only, in the
ambidextrous ``ideal'' generated by {$P,Q,K-1$}, but of course
this is not an ideal at all. It would be very nice if one had
a Gr\"obner basis algorithm for doing that. Nobuki Takayama made
considerable progress ({\it ``An approach to the zero recognition
problem by Buchberger's algorithm"}, 
J. Symbolic Computation {\bf 14}, 265-282 (1992).
 
\medskip
Let a discrete function $F(n,k)$ be annihilated by two operators
$P$ and $Q$, that are ``independent'' in some technical sense
(i.e. the form a holonomic ideal, see [Z1], [Ca]). Performing the
elimination process above (and the holonomicity guarantees that
we'll be successful), we get the operators $A,B,C$ and $S(N,n)$.
Now let
$$\displaylines{
G(n,k)=C(N,K,n,k) F(n,k),\cr
\noalign{\hbox{we have}}
S(N,n) F(n,k) = (K-1) G(n,k).\cr
\noalign{\hbox{It follows that}}
a(n):= \sum_{k} F(n,k) ,\cr
\noalign{\hbox{satisfies}}
S(N,n) a(n) \equiv 0 .\cr}
$$
 
Let's apply the elimination method to find a recurrence
operator annihilating $a(n)$, with
$$
F(n,k) := { {n} \choose {k} } 
{ {b} \choose {k} } =
{ {n!b!}  \over { k!^2 (n-k)! (b-k)! }} ,
$$
and thereby prove and discover the Vandermonde-Chu identity.
We have
$$\displaylines{
{{F(n+1,k)} \over {F(n,k)}}={{(n+1)} \over {(n-k+1)}} ,\cr
{{F(n,k+1)} \over {F(n,k)}}={{(n-k)(b-k)} \over {(k+1)^2 }}.\cr
}$$
Cross multiplying,
$$\eqalignno{
(n-k+1) F(n+1,k)  -(n+1)F(n,k)&  = 0
,\qquad&(i)\cr
(k+1)^2 F(n,k+1)  -  (n-k)(b-k) F(n,k)& = 0.&(ii)\cr
\noalign{\hbox{In operator notation:}}
[(n-k+1)N-(n+1)]F& =  0,\cr
[ (k+1)^2 K - (nb-bk-nk+ k^2 ) ] F& = 0.\cr}
$$
So $F$ is annihilated by the two operators $P$ and $Q$, where
$$
P=(n-k+1)N-(n+1)\, ; \quad Q= (k+1)^2 K - (nb-bk-nk+ k^2 ).
$$
We would like to find a good operator that annihilates $F$. By
{\it good} we mean ``independent of k'', modulo $(K-1)$ (where
the multiples of $(K-1)$ that we are allowed to throw out are
right multiples).
 
Let's first write $P$ and $Q$ in descending powers of $k$,
modulo $(K-1)$:
$$
P= (-N) k + (n+1)N-(n+1)\,;\quad
Q=(n+b)k-nb+(K-1) k^2 \,;
$$
and then eliminate $k$ modulo $(K-1)$.  However, we must be
careful to remember that left multiplying a general operator~$G$
by $(K-1) \hbox{\eightrm JUNK}$ does not yield, in general, 
$(K-1) \hbox{\eightrm JUNK}'$. In other words,
 
\bigskip
\noindent
{\it  Warning:}
$$
\hbox{\eightrm OPERATOR}(N,K,n,k) (K-1) 
(\hbox{\eightrm JUNK}) \neq (K-1) (\hbox{\eightrm JUNK}') .
$$
Left multiplying $P$ by $n+b+1$,  left multiplying $Q$ by $N$
and adding yields
$$\eqalign{
(n+b+1) P + N Q a&= (n+b+1) [-Nk+(n+1)N-(n+1)]\cr
&\qquad\qquad\qquad+{}
N[(n+b)k-nb+(K-1) k^2 ]\cr
&=(n+1)[(n+1)N-(n+b+1)]+(K-1)[N k^2 ] .\cr}
$$
So, in the above notation,
$$\displaylines{\hfill
S(N,n)=(n+1)[(n+1)N-(n+b+1)],\qquad 
\overline R  =  N k^2 .
\hfill\llap{(*)}\cr
\noalign{\hbox{It follows that}}
a(n):= \sum_{k} 
{{n} \choose {k} }
{{b} \choose {k} }\cr
\noalign{\hbox{satisfies}}
((n+1)N-(n+b+1)) a(n) \equiv 0 ,\cr}
$$
or, in everyday notation,
$$\displaylines{
(n+1)a(n+1)-(n+b+1)a(n) \equiv 0 ,\cr
\noalign{\hbox{i.e.,}}
a(n+1) = {{n+b+1} \over {n+1}} a(n) \,   \Rightarrow
\,
a(n)= { {(n+b)!} \over {n!}} C,\cr}
$$
for some  constant
 independent  of $ n$,
and plugging in $n=0$ yields that $1=a(0)=b!\,C$ and hence
$C=1/b!$. We have just discovered, and proved at the same time, the
Vandermonde-Chu identity.
 
\medskip
Note that once we have found the eliminated operator $S(N,n)$ 
and the corresponding $\bar R $ in $(*)$ above, we can present
the proof without mentioning how we obtained it. In this case
$\overline R = N k^2 $, so in the above notation
$$
G(n,k)  =   -\overline R  F(n,k)   =   - N k^2 F(n,k)  =  
{{-(n+1)!b!} \over {(k-1)!^2 (n-k+1)! (b-k)!} }     .
$$
So all we have to present are $S(N,n)$ and $G(n,k)$ above and
ask the readers to believe or prove for themselves the purely
routine assertion that
$$
S(N,n)F(n,k)=G(n,k+1)-G(n,k) .
$$

\smallskip
\noindent
{\it  Dixon's Identity by Elimination:}

We will now apply the elimination procedure to derive and prove
the celebrated Dixon identity of 1903. It states that
$$\displaylines{
\sum_{k} (-1)^k
{ {n+a} \choose {n+k} }
{ {n+b} \choose {b+k} }
{ {a+b} \choose {a+k} }
 = {{(n+a+b)!} \over {n!\,a!\,b!\,}}.\cr
\noalign{\hbox{Equivalently,}}
\quad
\sum_{k} {{(-1)^k }
\over
{(n+k)!\,(n-k)!\,(b+k)!\,(b-k)!\,(a+k)!\,(a-k)!\,}} \hfill\cr
\hfill{}
= {{(n+a+b)!\,} \over
{n!\,a!\,b!\,(n+a)!\,(n+b)!\,(a+b)!\,}} .\quad\cr}
$$
Calling the summand on the left $F(n,k)$, we have
$$\eqalign{
{{F(n+1,k)} \over {F(n,k)}}& = {{1} \over {(n+k+1)(n-k+1)}} ,\cr
{{F(n,k+1)} \over {F(n,k)}}& =
{{(-1)(n-k)(b-k)(a-k)} \over {(n+k+1)(b+k+1)(a+k+1)}} .\cr}
$$
It follows that $F(n,k)$ is annihilated by the operators
$$\displaylines{
P=N(n+k)(n-k)-1 ,\quad
Q=K(n+k)(a+k)(b+k)+(n-k)(a-k)(b-k)  .\cr
\noalign{\hbox{Rewrite $P$ and $Q$ in descending powers of $k$,
modulo $K-1$:}}
P  =   -N k^2 + (N n^2 -1 )  ,\cr
Q= 2(n+a+b) k^2 + 2nab + (K-1) ((n+k)(a+k)(b+k)).\cr}
$$
Now eliminate $k^2$, to get the following operator that
annihilates $F(n,k)$:
$$\displaylines{\noalign{\vskip-8pt}
\quad
2(n+a+b+1)P+NQ \hfill\cr
\kern2cm{}
=2(n+a+b+1)( N n^2 -1 ) + N(2nab)\hfill\cr
\kern3cm{}+
(K-1) ( N(n+k)(a+k)(b+k))  ,\hfill\cr
\noalign{\hbox{which equals}}
N[2n(n+a)(n+b)]-2(n+a+b+1)+ (K-1) (
N(n+k)(a+k)(b+k)).\cr}
$$
In the above notation we have found that the following operator
annihilates $a(n):= \sum_{k} F(n,k)$:
$$\displaylines{
\eqalign{
S(N,n)&=N[2n(n+a)(n+b)]-2(n+a+b+1)\cr
&=2(n+1)(n+a+1)(n+b+1)N
-2(n+a+b+1) . \cr}\cr
\noalign{\hbox{Also}}
\overline R  (N,K,n,k) = ( N(n+k)(a+k)(b+k)) ,\cr
\noalign{\hbox{and}}
\eqalign{
G(n,k) &=  - \overline R  F(n,k)\cr
&  = 
{{(-1)^{k-1} }
\over
{(n+k)!\,(n+1-k)!\,(b+k-1)!\,(b-k)!\,(a+k-1)!\,(a-k)!\,}}.
\cr}\cr}
$$
Once we have found $S(N,n)$ and $G(n,k)$ all we have to do is
present them and ask the readers to verify that
$$
S(N,n)F(n,k)=G(n,k+1)-G(n,k) .
$$

\smallskip
\noindent
{\it  Homework :}
 
1. Using the elimination method of this recitation find a
recurrence satisfied by
$$
a(n):= \sum_{k=0}^{n} 
{{ n-k} \choose {k} }  .
$$
(No credit for other methods!)
 
2. Find a recurrence satisfied by
$$
a(n):= \sum_{k} 
{ {n} \choose {k} } 
{ {n+k} \choose {k} }   .
$$
 
3(*). Using the method of this recitation, evaluate, if possible, 
the following sum :
$$
a(n):= \sum_{k} 
{{(a+k-1)!\,(b+k-1)!\,(c-a-b+n-k-1)!}
\over
{k!\,(n-k)!\,(c+k-1)!\,}}.
$$
If you succeeded you would have rediscovered and reproved the
Pfaff-Saalsch\"utz identity.
 
4.(**) Prove
$$
\sum_{  k_1 , k_2 }
{{(-1)^{k_1 + k_2 } (k_1+k_2)! }
\over
{ k_1 !^2 k_2 !^2 (n- k_1 )!\, (m- k_2 )!}}
 \, = \, 
c(n) \delta_{n,m}  .
$$
 
5. (100 F) Using elimination prove E3376 (AMM, March 1990):
$$
\sum_{i=0}^{ n} \sum_{j=0}^{n}
{ {i+j} \choose {j} }^2
{ {4n-2i-2j} \choose {2n-2j}}  = 
(2n+1) { {2n} \choose {n} }^2  .
$$
 
Note: This was partially solved by Peter Paule, see a forthcoming
paper joint with George Andrews [``J. Symbolic Computation" {\bf
16}, 147-153  (1993)].
Then it was completely solved by Peter Paule [``Solution of a
Seminaire Homework Example (28th SLC)'', RISC-Linz Report
Series No. 92-52, (1992).]
 
\bigskip
\centerline {\bf Recitation II. Gosper's Algorithm: A Decision
Procedure} 
\smallskip
\centerline{\bf for Indefinite Hypergeometric summation}
 
\medskip
Here I will describe and motivate Gosper's algorithm
[Proc. Nat. Acad. Sci. USA, 40-42 (1979)]. 
As will be explained in Recitation 3, Gosper's algorithm
for {\it indefinite summation} turned out to be even more important
for {\it definite summation} .

As we all know, a series 
$\sum\limits_{n} a(n)$ is called {\it geometric} if the ratio of
consecutive terms  are constant:
$$\displaylines{
{{a(n)} \over {a(n-1)}}  =   \hbox{\eightrm CONSTANT}.\cr
\noalign{\hbox{It is called {\it hypergeometric} if}}
{{a(n)} \over {a(n-1)}}   = 
\hbox{\eightrm RATIONAL FUNCTION OF }n.\cr}
$$ 

\goodbreak
\noindent
The sequence $\{ a(n) \} $ itself is called a {\it hypergeometric
sequence}, or, more often {\eightrm CLOSED FORM} (or \CF\ for
short). It is easy to see that every CF sequence can be
expressed as
$$
\hbox {\eightrm RATIONAL  FUNCTION} (n) z^n \,  \cdot 
{{\prod\limits_{i} ( a_i n + b_i ) !}
\over
{\prod\limits_{j} ( a'_j n + b'_j ) !}}  .
$$
Given a \CF, $a(n)$, Gosper asked, and brilliantly answered,
whether $S(n):= \sum\limits_{i=0}^{n} a(i)$ 
is also \CF, modulo a constant.
 
This is the discrete analog of Liouville's problem of ``integration
in finite form''. Since the discrete is much harder, and composition 
of discrete functions is badly behaved, we must be content with a
much narrower definition of \CF. The continuous counterpart of
what we call \CF\ would be functions $f(x)$ whose logarithmic
derivatives are rational functions, and hence functions of the
form
$$
\exp ( R_0 (x))
\prod_{i} R_i (x)^{\lambda_i } ,\quad R_i\hbox{ rational,}
$$
which is much narrower than Liouville's definition that allows
algebraic functions and compositions.

\medskip
\noindent
Going back to Gosper's problem, it can be phrased as follows.
 
{\bf Input:} \CF\ sequence $a(n)$.
 
{\bf Output:} \CF\ $S(n)$ such that
$S(n)  - S(n-1) = a(n) $ ,
or the statement ``does not exist''.
 
\medskip
Of course, {\it 
$S(n) \,  \CF $ implies $ a(n) \, \CF$.}

\medskip 
{\it  Proof:}
$$\eqalign{
{{a(n)} \over {a(n-1)}}  & ={{S(n)-S(n-1)} \over 
{S(n-1)-S(n-2)}} \cr
& = 
{{(S(n)/S(n-1) \, - \, 1)} \over
{(1 \, - \, S(n-2)/S(n-1))} }\cr 
&= \hbox{\eightrm  RATIONAL} \, (n).\qed\cr}
$$

You can find many ``good'' $a(n)$ by working backwards. Start with
a \CF\ $S(n)$, compute $a(n):= S(n) -  S(n-1)$,
and compare the forms of $S(n)$ and $a(n)$.
 
\medskip
\noindent
{\it  Example $1$}: $S(n)=n!$. Then
$$
a(n)=n!-(n-1)!=(n-1) (n-1)! .
$$
Note that $a(n)$ has two parts: $(n-1)$, the ``polynomial part'', 
which we call $p(n)$, and $(n-1)!$, the ``pure factorial part''.
Gosper's algorithm depends on such a decomposition. In
anticipation of Gosper's algorithm let us see  how $a(n)/a(n-1)$
looks like:
$$
{{a(n)} \over {a(n-1)}}
= {{n-1} \over {n-2}} \cdot {{(n-1)} \over {1}}
={{p(n)} \over {p(n-1)}} {{q(n)} \over {r(n)}}  .
$$
The $\displaystyle {p(n)\over p(n-1)}$ is there because of the
polynomial part, and the $\displaystyle {q(n)\over r(n)}$ is due
to the ``pure factorial part''. Anyway, with these names for the
parts of $\displaystyle {a(n)\over a(n-1)}$, we get that
$S(n)=n!$, in terms of $a(n)=(n-1)(n-1)!$, is 
$$
S(n)= {{a(n) n} \over {n-1}} =   {{a(n) q(n+1)} \over {p(n)}}. 
$$

\smallskip
\noindent
{\it  Example $2$:}  $S(n) \, =(n+3) n!$.
 
We have,
$a(n)=S(n)-S(n-1)=(n+3)n!-(n+2)(n-1)!=[(n+3)n-(n+2)](n-1)!=$
$( n^2 + 2n-2)(n-1)!$. Here the ``polynomial part'', $p(n)$, is
$ n^2 +2n-2$ and the ``pure factorial part'' is $(n-1)!$.
Now
$$
{{a(n)} \over {a(n-1)}} 
= {{p(n)} \over {p(n-1)}} {{(n-1)} \over {1}}.
$$
In anticipation of things to come, and as in the previous
example, let us call $(n-1)$ above $q(n)$ and $r(n)$, $1$. In
other words, if we write $a(n) = p(n) \overline a (n) $, where
$p(n)$ is the polynomial part and $\overline a (n)$ is the pure
factorial part, then
$$
{{q(n)} \over {r(n)}} 
:= {\overline a  (n) \over \overline a  (n-1)} .
$$
Recall that now we are working backwards, and that we already
know the answer $S(n)=(n+3)n!$. Let's see how it is expressible
in terms of $a(n)$ and its derived quantities, $p(n), q(n), r(n)$:
$$
S(n)=(n+3)n!=(n-1)!(n+3)n= \overline a  (n) (n+3) q(n+1) ,
$$
where $\overline a  (n)$ is the ``pure factorial part'',
$\displaystyle {a(n)\over p(n)}$. Thus it was possible to write
$S(n)= \overline a  (n) q(n+1) f(n)$, for some polynomial $f(n)$,
in this case of degree $1$. We will see that this is always
possible, and forms the essence of Gosper's algorithm.
 
\medskip
The above examples motivate the following way of ``guessing the
answer"
$$
S(n)= {{a(n)   q(n+1)} \over {p(n)}} f(n)  ,
$$
where $p(n)$ is the ``polynomial part'' of $a(n)$, obtained in the
decomposition
$a(n)= p(n) \overline a  (n) $,
of $a(n)$ as a product of polynomial part $p(n)$ and ``pure
factorial part'' $\overline a  (n)$, and $q(n)$ is the numerator
of $\displaystyle 
{{\overline a  (n)} \over {\overline a  (n-1)} }$.
In other words, $p(n), q(n), r(n)$ are the polynomials featuring
in the writing of \smash{$\displaystyle {a(n)\over a(n-1)}$} as
$$
{{a(n)} \over {a(n-1)}} 
= {{p(n)} \over {p(n-1)}} {{q(n)} \over {r(n)}}
$$
and $p(n)$ is maximal w.r.t.
$a(n)/a(n-1)$ being able to be written thus. It can be seen that
$q(n)$ and $r(n)$ satisfy
$$
\gcd ( q(n), r(n+j) ) =1\ {\rm for \  every\  integer}\ j 
\geq  0 . $$
If not, there exists a $j  \geq 0$ such that
$$
g(n):=\gcd( q(n), r(n+j) ) \neq 1 .
$$
Let
$$\displaylines{\noalign{\vskip -8pt}
q'(n) := {{q(n)} \over {g(n)}},
\   r'(n) := {{r(n)} \over {g(n-j)}},
\ p'(n):=p(n)g(n)g(n-1) \ldots g(n-j+1)  .\cr
\noalign{\hbox{Of course,}}
{{a(n)} \over {a(n-1)}} = {{p'(n)} \over {p'(n-1)}} {{q'(n)} \over
{r'(n)}}.\cr} 
$$
The above procedure gives an effective and efficient way to find
$p(n)$, $q(n)$, $r(n)$. Start with $p(n):=1$ (or rather with the
polynomial factor in front of $a(n)$), and get an initial
decomposition $$
{{a(n)} \over {a(n-1)}} 
= {{p(n)} \over {p(n-1)}} {{q(n)} \over {r(n)}}.
$$
Now check whether there exists a $j  \geq  0$ such that
$q(n)$ and $r(n+j)$ have a common factor. To find whether there exists
such a $j$, let
$$
R(j):={\rm Resultant}_n (q(n),r(n+j)),
$$
and find the non-negative integer roots of $R(j)=0$. In most
applications $q(n)$ and $r(n)$ come already factored:
$$
q(n)= \prod_{\alpha} (n- \alpha ),\qquad
r(n)= \prod_{\beta} (n- \beta ) .
$$
In this case it is easier to compute all 
the differences $\beta - \alpha$ and see if there is a
non-negative integer amongst them.

\medskip 
Sooner or later, we would arrive at a decomposition
$$
{{a(n)} \over {a(n-1)}} 
= {{p(n)} \over {p(n-1)}} {{q(n)} \over {r(n)}},
$$
with $\gcd( q(n), r(n+j))=1$ for every integer $j  \geq 0$.
Motivated by the above experimentation, we set (i.e. make
a {\it change of dependent variables}):
$$
S(n)={{ a(n)q(n+1)} \over {p(n)}} f(n).
$$
In the above, everything is known except $f(n)$. A priori, $f(n)$
is just another CF sequence, but the nice surprise is that:
 
\proclaim
Claim. The only way that $S(n)$ is \CF\ is for $f(n)$ to be a
rational function.
\endproclaim

{\it  Proof}:
$$\eqalign{
f(n)&= {{p(n)S(n)} \over {q(n+1)a(n)}}\cr
&  = 
{{p(n)S(n)} \over {q(n+1)(S(n)-S(n-1))}} \cr
&= 
{{p(n)} \over
{q(n+1)(1-S(n-1)/S(n))}} ,\cr}
$$
and thus must be a rational function, if $S(n)$ is \CF.\halmos
 
\medskip
What does $S(n)-S(n-1)=a(n)$
say about $f(n)$?. It is easily seen that the equation
for $f(n)$ is
$$
q(n+1) f(n) - r(n)f(n-1)  = p(n) .
\eqno(*)
$$
This is the ``{\eightrm FUNCTIONAL EQUATION FOR} $f(n)$''.
 
\proclaim
Surprise. The only way that $f(n)$ can be a rational function is
for it to be a polynomial.
\endproclaim
 
{\it  Proof:} A starred homework exercise.
 
{\it  Hint:} Suppose $f(n)=c(n)/d(n)$, $d(n)  \neq  1$. Let $j$ be
the largest integer such that $gcd(d(n),d(n+j) = g(n)  \neq  1$,
and arrive at a contradiction from the functional equation, the
assumption on $q(n),r(n)$ and the maximality of $j$.
 
How to solve the functional equation? We need an upper
bound for the degree of $f(n)$. Equating degrees, we get
$$\displaylines{
\deg f+ \max( \deg q, \deg r)  =  \deg p .\cr
\noalign{\hbox{So unless there is some fluke,}}
L:=\deg f= \deg p-\max(\deg q,\deg r)  .\cr}
$$
The fluke happens when the two leading coefficients of
$q(n+1)$ and $r(n)$ are such that it is possible for a higher
degree polynomial $f(n)$ to exist, which will make the leading
coefficient of the left side vanish, and hence make it still
possible for the degree of $f(n)$ to be higher. This must be
checked, and then one has to take a larger $L$. All this is
described in Gosper's paper.  For pedagogical reasons, we won't
worry about it here. However, as pointed out by Petr Lisonek,
Peter Paule, and Volker Strehl, this case comes up pretty often,
especially in the context of the fast algorithm. See their paper:
{\it ``Improvements of the degree settings in Gosper's 
algorithm''}, JSC {\bf 16}, 243-253 (1993).
 
Having found an upper bound for the degree $L$ of $f(n)$, we set
$$
f(n)= \sum_{i=0}^{L} f_i n^i  ,
$$
plug into $(*)$, compare coefficients and solve the resulting
system of linear equations.

\medskip
\noindent 
{\it  Example 1:} Find out whether
$\sum\limits_{n} (n-1) (n-1)!$ has closed form.
 
\medskip 
{\it  Solution:} Here $a(n)=(n-1)(n-1)!$.
%{\bf Step 1:}
Step 1 is:
$$
{{a(n)} \over {a(n-1)}} 
= {{n-1} \over {n-2}} \cdot {{n-1} \over {1}},
$$
so initially, $p(n)=n-1$, $q(n)=n-1$, and $r(n)=1$. Obviously,
$\gcd(q(n),r(n+j))=1$, for every $j  \geq  0$, so these values for
$p(n)$, $q(n)$, $r(n)$ are the final ones. The functional equation
reads
$$
nf(n) - f(n-1) = n-1  ,
$$
$L:=\deg f=1-1=0$, so $f= f_0$. Plugging this into the
functional equation we get
$$
n f_0 - f_0 = n-1.
$$
Equating coefficients of $n$ and $n^0$, we get the
two equations
$$
f_0 =1 ,\quad - f_0 = -1  .
$$
The solution is $f_0 =1$, so $f(n)=1$, and thus
$$
S(n)= {{a(n) q(n+1) f(n)} \over {p(n)}} =
{{(n-1)(n-1)! \cdot n \cdot 1 } \over {(n-1)}} = n!  .
$$
Checking we see that indeed $n!-(n-1)! = (n-1)(n-1)!$.
 
\medskip\noindent
{\it  Example 2:} Is the sum 
$S(n):= \sum\limits_{i=0}^{n} i!$ expressible in closed form?

\medskip
{\it Solution:}
Here $a(n)=n!$, so $\displaystyle 
{{a(n)} \over {a(n-1)}} = n $.
Here $p(n)=1$, $q(n)=n$, and $r(n)=1$. The functional equation is
$$
(n+1)f(n)-f(n-1)=1 .
$$
This is impossible since the degree of $f$ should be $-1$.

\medskip \noindent
{\it  Homework:}
 
1. Is the  sum $\displaystyle \sum_{m=0}^{n} {{(2m)!} \over
{m!\,(m+1)!}}$ expressible in closed form? 
 
\noindent
(Ans.: No.)

2. (Amer. Math. Monthly, Nov. 1989, problem E3352). 
Prove
$$
\sum_{n=0}^{\infty}
{{1} \over 
{n!\,( n^4 + n^2 + 1 )}}  \, = \,  {{e} \over {2}}.$$
 
3. Can the harmonic numbers
$
H_n =\displaystyle  \sum_{i=1}^{n} {{1} \over {i}}$
be expressed in closed form? 
(You are supposed to use Gosper's algorithm, but it is possible to
prove this using asymptotics, as shown by Gilbert Labelle.)
 
4. We all know that, for any fixed $A$,
$\displaystyle 
\sum_{k=1}^{A} {{A} \choose {k}} = 2^A $.
Is there a closed form expression in $n$, for the partial sums
of the binomial coefficients
$S(n):=\displaystyle  
\sum_{k=0}^{n} {{A} \choose {k}} $?

5. Find, if possible,
$\displaystyle 
\sum_{n=0}^{m} {{(4n-3)(2n-2)!} \over {(n-1)!}} $.
 
\vfill\eject
\centerline
{\bf Recitation III:  From Indefinite Hypergeometric Summation}
\smallskip
\centerline{\bf To Definite Hypergeometric summation and WZ
Pairs}
 
\medskip
Gosper's algorithm for
{\it indefinite} summation is the basis for my algorithm for 
{\it definite} summation, but not in the obvious way! Most definite
identities
$$
\sum_{k = - \infty}^{\infty} F(n,k)  
 = \hbox{\eightrm NICE}
(n) 
\quad {\rm have}\quad
\sum_{k = - \infty}^{m} F(n,k)  =\hbox{\eightrm  UGLY} (n,m).$$
If
\smash{$\displaystyle 
\sum_{ k = - \infty}^{m} F(n,k)  = \hbox{\eightrm NICE} (n,m)$},
then Gosper's method can be used to find 
$\hbox{\eightrm NICE} (n,m)$, and
$\displaystyle 
\sum_{k = - \infty}^{\infty } F(n,k)  
= \hbox{\eightrm NICE} (n, \infty ) = \hbox{\eightrm NICE}(n)$.
Whenever that is the case the definite identity is {\it trivial}.
To take a metaphor from calculus,
$\displaystyle 
\int_{- \infty }^{ \infty } e^{ - x^2 } dx  =  \sqrt{\pi }$
is deep since the corresponding indefinite integral
$\displaystyle 
\int_{ - \infty }^{ x } e^{ - t^2 }\, dt$
is not expressible in closed form, while
$\displaystyle 
\int_{ - \infty }^{ \infty } x e^{ - x^2 }\, dx = {{1} \over
{2}} $
is shallow, since the integrand has an antiderivative that is
expressible in closed form.

\medskip
My fast algorithm starts with a definite sum
$a(n):= \sum\limits_{k} F(n,k)$
and finds a homogeneous
linear recurrence equation with polynomial coefficients
satisfied by $a(n)$. If the recurrence is first order, then
$a(n)$ can be easily expressed explicitly, otherwise we must be
content with the recurrence. The algorithm is not guaranteed to find
the minimal recurrence, although it usually does. 
Marko Petkovsek has recently come up with
a beautiful algorithm that decides when a linear recurrence has
closed form solutions. The combination of my fast algorithm and
Petkovsek's algorithm 
[{\it Hypergeometric solutions of linear recurrence
equations with polynomial coefficients}, J. Symbolic Computation
{\bf 14}, 243-264 (1992)] completely solves the problem of
deciding when a {\it definite} hypergeometric sum can be
expressed in closed form.
 
\medskip
Let's first consider the special case of sums 
$a(n) = \sum\limits_{k} F(n,k)$, for which $a(n)$ satisfies a first
order recurrence, so that one has an ``identity''. For that
important special case, Herb Wilf made a brilliant observation
that at first only seemed to be a minor simplification, and like
all great discoveries, seems obvious by hindsight, but it led to
the conceptual breakthrough of {\it WZ pairs} [WZ1-2], and {\it
WZ forms} [Z4]. 
 
\vfill\eject

\proclaim Wilf's brilliant idea. Instead of trying to prove
\smash{$\sum\limits_{k} F(n,k)  = \hbox{\eightrm NICE} (n)$} try
to prove $\displaystyle\sum_{k}  {{F(n,k)}  \over
{\hbox{\eightrm NICE}(n)}}=1$.
\endproclaim

Renaming the summand on the 
left side of the above $F(n,k)$, we are left with
the task of proving, for given Closed Forms $F(n,k)$, identities
of the form $
\sum\limits_{k} F(n,k)=1$.
 
Let us call the left side $a(n)$. We have to prove that $a(n)
\equiv 1$. It is always trivial to check, in any given instance,
that $a(0) = 1$. The assertion that $a(n) \equiv 1$ would then
follow by induction if we can show  that 
$$a(n+1)  -a(n)  \equiv  0,\quad{\rm i.e.,}\quad
\sum_{k} (F(n+1,k)-F(n,k)) =0.
$$
 
\proclaim Big Surprise {\rm (Gosper's Missed Opportunity)}.  
Although $\sum\limits_{k} F(n,k)$ is (usually) not indefinitely
summable, in the vast majority of cases,
$$\displaylines{%\noalign{\vskip -8pt}
\sum_{k}
(F(n+1,k)-F(n,k))\cr
\noalign{\hbox{is!, i.e., there exists a closed form $G(n,k)$
such that}}
\hfill
F(n+1,k)-F(n,k) =G(n,k+1)-G(n,k)  .
\hfill\llap{\rm (WZ)}\cr}
$$
\endproclaim

The pair $(F,G)$ is called a {\it WZ pair}. To prove
$\sum\limits_{k} F(n,k) \equiv 1$, all we have to do is present
the ``certificate''  $G$, and the reader can then check that
$(F,G)$ is a WZ pair. The proof then follows upon summing (WZ) 
with respect to $k$.
 
Thanks to Gosper's algorithm, we can always find the $G(n,k)$
whenever it exists, and we know that its form is
$$\displaylines{%\noalign{\vskip -8pt}
\quad
\hbox{\eightrm RATIONAL}(n,k)[F(n+1,k)-F(n,k)] \hfill\cr
\kern 2.5cm{}=  \hbox{\eightrm
RATIONAL}(n,k)[F(n+1,k)/F(n,k)-1]F(n,k)\hfill\cr 
\kern 2.5cm{}=R(n,k) F(n,k)\
{\rm (say)}.\hfill\cr} 
$$
Hence it is enough to give the \hbox{\eightrm RATIONAL}
function $R(n,k)$.
 
\medskip\noindent
{\it  Example}: $\displaystyle \sum_{k} 
{ {n}  \choose {k}}  =2^n$.
Here $F(n,k) =\displaystyle {  n!\over k!\,(n-k)! \,2^n} $
and\hfil\break $F(n+1,k)-F(n,k) =\displaystyle  
{{ -(n-2k+1)n!} \over { 2^{n+1} k!\, (n-k+1)!}}$.
Using Gosper's algorithm, we find that
the antidifference of this w.r.t. $k$ is
$$
{\rm nusum} (\%,k) =  {{-n!} \over { 2^{n+1} (n-k+1)!(k-1)!}}
.$$
So
$$
G(n,k) = {{ - 1} \over { 2^{n+1} } }
{ {n} \choose {k-1} } 
 .
$$
 
\goodbreak
\proclaim Bonus. Buy one identity and get one identity free.
\endproclaim

Summing (WZ)  w.r.t. $n$, we get
$$\displaylines{
0  =  \sum_{n} (  F(n+1,k)-F(n,k)  )  =  
\sum_{n} (G(n,k+1)-G(n,k)),\cr
\noalign{\hbox{and hence}}
\sum_{n} G(n,k)  =  C ,\cr}
$$
where $C$ is a constant independent of $k$, which can be easily
evaluated by plugging in $k = 0$. This is called the {\it dual
identity}.
 
In practice the above procedure will yield $C =  \infty$, i.e., the
sum diverges. However one can get  new non-trivial identities in
two different ways. The first one is by summing, not from $n= -
\infty$ to $n= \infty$, but rather from $n=0$ to $n= \infty$.
When we do that, we get
$$
\sum_{n=0}^{\infty} G(n,k)   
=   \sum_{j  \leq k-1} ( f_j - F(0,j))  .
$$
Here $f_j$ is defined by
$\displaystyle 
f_j := \lim_{n \rightarrow \infty }  F(n,j)$,
which is usually a triviality to compute. I refer the reader to
[WZ1], cited at the beginning of these notes, for several
interesting examples. In addition, the identities in the
very last section of Bailey's book {\it Generalized
Hypergeometric Series}, that seemed hitherto mysterious and
artificial, all emerge as companion identities of well known ones.
 
A second way of obtaining a companion identity is by
introducing ``shadows''. This has the advantage that one still gets
standard identities in which the right hand side has closed form.
 
\medskip 
\proclaim Shadow. The operation of shadowing is like discrete
``analytic continuation''. The expression $n!$ is meaningless for
$n$ negative, or if you wish has a singularity there. 
\endproclaim

But what makes $n!$ what it is? The defining property is that
$a(n):=n!$ satisfies the recurrence equation
$a(n) = na(n-1)$, with the initial condition $a(0) = 1$. If we try to use
it to define the value of $a(n)$ at $n = -1$, by plugging $n = 0$, we get
$a(0) = 0a(-1)$. So there is no function $a(n)$ that is defined for
all integers $n$ and that satisfies $a(n) = na(n-1)$. But what is so great
about the positive integers? We can ask that $a(n) = na(n-1)$ holds
for negative integers! We get $\overline a  (n) =\displaystyle  
{{(-1)^n} \over {(-n-1)!}}$.
 
 We call $\overline a  (n)$ above the {\it shadow} of $n!$. It
satisfies the same recurrence as that of $n!$, but is defined for
the set of negative integers rather than positive integers.
 
More generally, the {\it shadow} of a factorial of a linear
expression: $(an+bk+c)!$, with $a,b$ integers and $c$ any indeterminate,
is defined by
$$
(an+bk+c)!   \rightarrow   {{(-1)^{an+bk+c}} \over
{(-an-bk-c-1)!}}  .
$$
The shadow of $(an+bk+c)!$ satisfies the same linear
recurrence equations with polynomial coefficients as
$(an+bk+c)!$ since $F(n+1,k)/F(n,k)$ and $F(n,k+1)/F(n,k)$
give the same \hbox{\eightrm RATIONAL} functions respectively
for both $F(n,k) = (an+bk+c)!$ and $F(n,k) =  (-1)^{an+bk+c}
/(-an-bk-c-1)!$.  Thus everything that is true for one, as far
as elimination and Gosper's algorithm are concerned, is also
true for the other, and for the purposes of the present theory,
they are completely equivalent.  The only difference is in their
domain of definition, and when they vanish. 
 
Finally if one has $F(n,k)$ equal to a power times a quotient of products
of such linear terms, one can apply the shadow treatment to any
number of the terms $(an+bk+c)!$ that appear on either the numerator
or denominator, getting $2^{\hbox{\eightrm \#  of  such 
terms}}$ possibilities for equivalent $F(n,k)$. So if one has a
sum $\sum\limits_{k} F(n,k)$
which diverges for $n$, one can always find an equivalent
$F(n,k)$ for which the sum converges for a ``half discrete line''
in $n$.
 
In practice, the default shadowing of such a summand $F(n,k)$
would be obtained by shadowing each term $(an+bk+c)!$ for
which $a+b  \neq  0$ and leaving all terms of the form
$(an-ak+c)!$ alone.
 
Recall the WZ pair that arose above, when we proved that the sum
of the binomial coefficients $n!/(k!(n-k)!)$ was $2^n$:
$$
(F,G):=
\biggl( {{1} \over { 2^n }} { {n} \choose {k} } \ , \  
-{{1} \over { 2^{n+1} }} { {n} \choose {k-1} } \biggr).
$$
The dual sum $\sum\limits_{ n} G(n,k)$ diverges for every $k$. To
make it meaningful, consider the shadow WS {\it pair}:
$$
( \overline F  , \overline G ):=
\biggl (  {{(-1)^{n+k}} \over { 2^n }} {{ -k-1} \choose
{-n-1} }  \, , \, 
 {{(-1)^{n+k}} \over { 2^{n+1} } }
{ {-k} \choose {-n-1} }  \biggr )  .
$$
Now $\overline G (n,k)$ has compact support w.r.t. $n$ for all
negative $k$, and we deduce
$$
\sum_{ n}
{{(-1)^{n+k}} \over { 2^{n+1} }} { {-k} \choose {-n-1} } 
  =  C \quad\hbox{for  each  negative } k.
$$
Making the transformation $ k \leftarrow  -k$, $n \leftarrow 
-n-1$, we get
$$
\sum_{n} (-2)^n
{ {k} \choose {n} } = (-1)^k .
$$
So it turned out that the dual of $ (1+1)^n  =  2^n$ is
$(1-2)^k  =  (-1)^k$.

\goodbreak

\medskip\noindent
{\it  Exercise:} Find the dual identity to the binomial theorem
$$\sum_{k} { {n} \choose {k} } x^k  =  (1+x)^n.$$

\smallskip
Many identities have free parameters. By specializing we get
``new'' identities, that are trivially implied by the original, more
general identities. 
Now comes an {\it important empirical observation}:
 
\proclaim Observation. The dual of a specialization is not, in
general, a specialization of the dual.
\endproclaim 
 
It follows that one can crank out lots of brand new identities,
complete with proofs, that a priori
are highly non-trivial, by iterating specialization and dualizing.
 
\medskip \noindent
{\it  Example} 
{\eightrm (SPECIALIZE AND DUALIZE)*}:
 The general Vandermonde identity is
$$
\sum_{k} {{n} \choose {k}}
{{a} \choose {k}}=
{{n+a} \choose {a}}.
$$
Its dual identity is just another rendition of same, with
changed parameters. But now specialize $n = a$:
$$\displaylines{
\sum_{k} {{n} \choose {k}}^2
={{2n} \choose {n}} .\cr
\noalign{\hbox{The dual of this is (check!)}}
\sum_{k} (3k-2n) {{n} \choose {k}}^2 
{{2k} \choose {k}}=0.\cr}
$$
 
This is {\eightrm A BRAND NEW IDENTITY}, unknown to Askey. It has
a $q$-analog derived from the $q$-version of WZ, that was
unknown to Andrews, and even whose limiting case was brand new,
and it took George Andrews three densely packed pages, using
five different identities, to prove.

\medskip 
\noindent
{\eightrm WHAT IS THE SECRET BEHIND THE WZ MIRACLE?}
 
\smallskip
If $F(n,k)$ is Closed Form, it is holonomic. Indeed, we have that
$$
{{F(n+1,k)} \over {F(n,k)}}   =   {{A(n,k)} \over {B(n,k)}} 
,\quad
{{F(n,k+1)} \over {F(n,k)}}   =   {{C(n,k)} \over {D(n,k)}}      ,
$$
for some polynomials $A$, $B$, $C$, $D$. (Of course they must
satisfy the obvious compatibility condition,) so, introducing the
operators
$$
P:=B(n,k)N  -  A(n,k),\quad Q:=D(n,k)K  -  C(n,k)    ,
$$
we see that $F(n,k)$ is annihilated by both $P$ and $Q$. 
By the first lecture, we know
that there exist operators $X(N,K,n,k)$ and $Y(N,K,n,k)$  and
$Z(N,K,n,k)$ such that
$$\displaylines{\quad
S(N,n):=X(N,K,n,k)P(N,K,n,k)+Y(N,K,n,k)Q(N,K,n,k)\hfill\cr
\hfill{}+(K-1)C(N,K,n,k)\quad\cr}
$$
is independent of $K$ and $k$.
Calling $
G(n,k):=C(N,K,n,k)F(n,k)$,
we get that
$S(N,n)F(n,k) = (K-1)G(n,k)$,
or in everyday notation,
$$S(N,n)F(n,k) = G(n,k+1)-G(n,k).$$

\proclaim Important Observation. If $F(n,k)$ is Closed Form, so
is $G(n,k)$.
\endproclaim

{\it  Proof}:
$$\eqalign{
N^i K^j F(n,k) & =  F(n+i,k+j)  = 
{{F(n+i,k+j)} \over {F(n,k)}}
 \cdot F(n,k)  \cr
&=  [ \hbox{\eightrm RATIONAL}(n,k) ] F(n,k)  .\cr} 
$$
Since, for any operator $C(N,K,n,k)$, $C(N,K,n,k)F(n,k)$ is
a linear combination, with coefficients that are polynomials
in $n$ and $k$, of terms as above, it follows that
$$
C(N,K,n,k) F(n,k)  =  
\hbox{\eightrm RATIONAL}(n,k) F(n,k) .\qed
$$

Going back to proving identities of the form
$\sum\limits_{k} F(n,k)=1$,
we want to prove that
$a(n):= \sum\limits_{k} F(n,k)$
satisfies the recurrence 
 $(N-1)a(n) =0$.

The elimination algorithm gives {\it a} recurrence
$S(N,n)a(n)  \equiv  0$,
that came from
$$
S(N,n)F(n,k)  \equiv   G(n,k+1)-G(n,k),
$$
for some closed form $G(n,k)$ (that is a multiple of
$F(n,k)$ by a \hbox{\eightrm RATIONAL} function). 
Let the order of $S(N,n)$ be {\eightrm ORDER}. To complete the
proof that $a(n) \equiv 1$, all we have is to check that this  is
true for $n = 0, \ldots , \hbox{\eightrm ORDER} -1$, and then
check that $S(N,n) 1  \equiv  0$. Equivalently, we have to see
whether $S(N,n)$ is a left multiple of $N-1$.
 
The WZ miracle takes place exactly when
the elimination algorithm actually gives us  $S(N,n) = N-1$, and
not a left multiple of it. It turns out that in the vast majority
of cases we are lucky, and for those cases, it  suffices to have
the WZ theory, and not the more general theory behind it.
However, 
 
(i) Sometimes we are not lucky, and $S(N,n)$ is not first order
 
(ii) What if we don't know the answer? In WZ theory, you should
know or guess, the answer.
 
(iii) What if the sum doesn't evaluate in closed form. The general
holonomic machinery promises us that the sum satisfies a linear
recurrence equation with polynomial coefficients, that should be
possible to find by elimination, using the method of Recitation I.
However, elimination is very slow.
 
\medskip\noindent
The question is:

{\eightrm  
IS THERE A FAST ALGORITHM
FOR FINDING THE RECURRENCE $S(N,n)$ AND THE ACCOMPANYING 
``CERTIFICATE"}
$G(n,k)$?
 
\medskip\noindent
The answer is: {\eightrm YES}.
 
\medskip
A simplistic way would be to ``guess'' empirically the recurrence
$$S(N,n) a(n)  = 0$$ 
satisfied by $a(n)$ and then use Gosper's
algorithm, w.r.t. $k$ to find a closed form $G(n,k)$ such that
$$
S(N,n)F(n,k) = G(n,k+1)-G(n,k)  .
$$
However, this has two drawbacks. One is practical: we don't know
what the degrees of the coefficients of $S(N,n)$ are going to be, and
we have to keep trying bigger and bigger degrees. 
The other is philosophical: this is empirical guessing. Finally, we
are not guaranteed that it is going to work. We do know, for sure
that there exists an operator $S(N,n)$ s.t. $S(N,n)F(n,k) =
G(n,k+1)-G(n,k)$, for some closed form $G(n,k)$ that is a
multiple of $F(n,k)$ by a \hbox{\eightrm RATIONAL} function. This
implies that $S(N,n) a(n) \equiv 0$. 

But the converse is not true: $a(n)$ may satisfy
a lower order recurrence, $S_1 (N,n) a(n) \equiv 0$. This recurrence
will be found empirically, but Gosper's algorithm will fail when
we try to find $nusum ( S_1 (N,n) F(n,k),k)$. To conclusively
and rigorously prove that $S_1 (N,n) a(n) \equiv 0$, we 
(or rather our computers) ``divide'' $S(N,n)$ by $S_1 (N,n)$:
$S(N,n) =  T(N,n) S_1 (N,n)$, and make sure that there is
no remainder. Since the elimination algorithm guarantees that
$S(N,n) a(n) \equiv 0$, we know that $T(N,n) [ S_1 (N,n) a(n) ] =
0$, and hence $ S_1 (N,n) a(n) \equiv 0$, provided it is true for
the first few values of $n$, which we already know is true, since
we found $S_1 (N,n)$ empirically at the first place.
 
 So what we really want is a {\eightrm FAST} algorithm for finding
an operator $S(N,n)$ and an accompanying closed form function
$G(n,k)$ such that
$$
S(N,n)F(n,k)  =  G(n,k+1)-G(n,k) .
\eqno(*)
$$
 
Let's suppose that we already know $S(N,n)$ by other means, but
still have to find $G(n,k)$. Then, it can be found using Gosper's
algorithm! This follows from the fact that $S(N,n) F(n,k)$ is
closed form itself, as shown above, and hence Gosper's algorithm
with respect to $k$ would produce the closed form 
anti-difference $G(n,k)$ whenever it exists, (and it does exist
thanks to the assumption.)
 
The problem is that we {\it don't} know $S(N,n)$ beforehand. We have to
find {\it both} $S(N,n)$ and $G(n,k)$ at the same time, 
from {\it scratch}, starting from the input $F(n,k)$. 
The pleasant surprise is:
 
\medskip\noindent
{\eightrm GOSPER'S ALGORITHM CAN BE EXTENDED TO MANUFACTURE
BOTH\hfil\break $G(n,k)$ AND $S(N,n)$ AT THE SAME TIME!}
 
\smallskip
What we do is a little like Lagrange multipliers. We first ``guess"
the order $I$ of the recurrence, and write $S(N,n)$ in generic form
$$
S(N,n):= \sum_{ i=0}^{I} s_i (n) N^i ,
$$
where the coefficients $s_i (n)$, which are polynomials in $n$,
have to be determined. In practice there is no ``guessing'' at all, since
we start with $I=0$ and do-loop our way up until we are successful. The
general holonomic theory and elimination procedure of recitation 1 
guarantees us success eventually. Furthermore, it's possible to
give a priori upper bound for $I$.
 
We now work with the generic $s_i (n)$ as though we knew what
they were, and form
$$\eqalign{
H(n,k):=S(N,n)F(n,k) & =  
\sum_{i=0}^{I} s_i (n) F(n+i,k) \cr
& =\Bigl[ \sum_{ i=0}^{I} s_i (n) {F(n+i,k)\over F(n,k)}\Bigr]
F(n,k)  .\cr} $$
The quantity in square brackets is a certain {\it
\hbox{\eightrm RATIONAL} function}, whose numerator is a
{\eightrm LINEAR EXPRESSION IN THE} $s_0 (n) , \ldots , s_I (n)$.
 
Note that when we do Gosper w.r.t. $k$, $n$ is a mere auxiliary
parameter and all the calculations are done in the field of \hbox{\eightrm RATIONAL}
functions in $n$.
So, let's do Gosper w.r.t. $k$ and let's take a look
at the functional equation for the polynomial $f(k)$ that determines
the closed form anti-difference of $H(n,k)$ above:
 $$
p(k)  =  q(k+1)f(k)  -  r(k)f(k-1) .
$$
Recall that we find $f(k)$ by expressing it in generic form
$f(k) =  f_0 + f_1 k + \cdots + f_L k^L$, plugging in the
functional equation and comparing coefficients of respective powers of
$k$. This gives a system of linear equations in 
$f_0 , f_1 , \ldots , f_L$, where the right sides involve
certain expressions in the $s_i (n)$. In fact we are asking
for $s_i (n)$ that will make the equations solvable. The
miracle is that the $s_i (n)$ occur linearly. So what we
have, in fact, is a linear system of equations with unknowns
$f_0 , \ldots , f_L $ and $s_0 ,\ldots , s_I$. If the system is not
solvable, then it means that there does not exist any recurrence
$S(N,n)$, of order $I$, such that (*) is true, and we must try
again, replacing $I$ by $I+1$. The general proof above
guarantees that we are going to be successful eventually.
 
 \medskip\noindent
{\it Example:} Let's find a recurrence for
$a(n):=\displaystyle  \sum_{k} {{1} \over {k!(n-k)!}}$.

\smallskip 
Here $F(n,k) =  1/(k!(n-k)!)$. We first try $I = 0$ and fail
(this means that the corresponding indefinite sum does not exist
in closed form), we then try $I = 1$ and set
$S(N,n) =  s_0 (n) + s_1 (n) N$.
Now
$$\eqalign{
H(n,k) =  S(N,n) F(n,k)  &=  s_0 (n) F(n,k) + s_1 (n) F(n+1,k)\cr
& ={{s_0 (n)} \over {k!(n-k)!}} + {{s_1 (n)} \over
{k!(n+1-k)!}} \cr
&={{ s_0 (n) (n+1-k) + s_1 (n) } \over {k! (n+1-k)!}}.\cr}
$$ 
We now do Gosper's algorithm, as described in recitation II,
``pretending" that we know  what $s_0 (n)$ and $s_1 (n)$ are. 
Using the notation of the last recitation (which coincides with
Gosper's notation), we have initially
$$\eqalign{
p(k)&  =   (n+1-k) s_0 (n) + s_1 (n)   ,  \cr  
{{q(k)} \over {r(k)}}&   =   
{{((k-1)!(n-k+2)!} \over {k!(n+1-k)!}}  =  
{{(n+2-k)} \over {k}}      ,\cr}
$$
so initially $q(k) = n+2-k$ and $r(k) = k$. Now we must make
sure that
$$
\gcd(q(k),r(k+j)) = 1\ \hbox{for  every  integer } j  \geq  0.
$$
Since $n+2-k$ and $k+j$ never have a common factor, this is
certainly true, so the final $p(k)$, $q(k)$, and $r(k)$ are given
by
$$
p(k)  =   (n+1-k) s_0 (n) + s_1 (n) ,\quad
q(k)=n+2-k ,\quad r(k)=k    .
$$
Substituting this into Gosper's functional equation:
$$\displaylines{
q(k+1)f(k)  -  r(k)f(k-1)  =  p(k)    ,\cr
\noalign{\hbox{we get, in this case,}}
 (n-k+1)f(k)  -  kf(k-1)  =   (n+1-k) s_0 (n) + s_1 (n). \cr}$$
\goodbreak

\noindent
The degree, in $k$, of the right side, is $1$, while the degree
of the left side is $1+ \deg _k f$. So $\deg_k f  =  1-1 = 0$.
We thus set $f(k) := f_0$ above, and get
$$
 (n-k+1) f_0   -  k f_0   =   (n+1-k) s_0 (n) + s_1 (n) .
$$
Comparing coefficients of $k^1$ and $k^0$ respectively yields
$2$ homogeneous linear equations for the three unknowns
$f_0$, $s_0$ and $s_1$:
 $$\eqalignno{
(n+1) f_0 & =  (n+1) s_0 (n) + s_1 (n),
&(i)\cr
 - 2 f_0  &=  - s_0 (n) .&(ii)\cr}
$$
Normalizing $f_0  =  f = 1$, we get the solution:
$s_0 = 2$, $s_1 = -(n+1)$, and hence
$$
f(k)  =  1,\quad  S(N,n) =  s_0 + s_ 1 N  =  2  -  (n+1)N  .
$$
Implementing $f(k)$, we get that $G(n,k)$ of (*) is given by
$$\eqalign{
G(n,k)& =  {{H(n,k)} \over {p(k)}}  q(k+1) f(k)\cr
&  =  
{{1} \over {k!(n-k+1)!}} \cdot
(n-k+1) \cdot 1     =  {{1} \over {k!(n-k)!}} .\cr}
$$
We have just found the recurrence satisfied by
$a(n) := \sum\limits_{k} \displaystyle {{1} \over {k!(n-k)!}}$.
 
It is $(2-(n+1)N)a(n) = 0$, i.e. $2a(n)-(n+1)a(n+1) = 0$, so
$a(n+1) =  (2/(n+1))a(n)$ which implies the closed form answer
$a(n) =  2^n / n!$. The proof consists in presenting the ``proof
certificate'' $G(n,k) =  1/(k!(n-k)!)$, 
and urging the readers to verify,
or believe that
$$
 2F(n,k) - (n+1) F(n+1,k)  =  G(n,k) - G(n,k-1 .
$$
The proof then follows by summing w.r.t. $k$.
 
 \bigskip
\noindent
{\it Homework:}
Using the algorithm of this recitation, find  recurrences for
the following binomial coefficients sums:
$$\eqalignno{
&\sum_{k}^{n} {{n-k} \choose {k}},
&(1)\cr
&\sum_{k} {{n} \choose {k}}^2,
&(2)\cr
&\sum_{k}{{n} \choose {k}}
{{n+k} \choose {k}}.
&(3)\cr}
$$
 
\goodbreak
\medskip
 
\centerline{\bf Postscript}
 
\medskip
Everything here has been $q$-ified. There is a $q$-analog of
Gosper's algorithm, and of its extension described above, that
would have appeared in my paper ``The method of creative
telescoping for $q$-series'', that became unnecessary because of
Tom Koornwinder's brilliant paper [K]. There is also a continuous
analog that appeared in my paper with Almkvist [AZ] cited at the
beginning of these notes. 

The next step would be to find a
{\eightrm FAST} algorithm for multisums. It follows from the
general holonomic theory that whenever $F(n, k_1 , \ldots , k_r )$ is
closed form, there exists an operator $S(N,n)$, and
closed form 
$$\displaylines{
G_1 (n, k_1 , \ldots , k_r ), \ldots,
G_r (n, k_1 , \ldots , k_r ),\cr 
\noalign{\hbox{such that}}
S(N,n) F(n, k_1 , \ldots , k_r ) 
 =  [G_1 (n, k_1 , \ldots , k_r ) -
G_1 (n, k_1 -1 , \ldots , k_r )] \hfill\cr
\hfill{}+ \cdots +
[G_r (n, k_1 , \ldots , k_r ) -
G_1 (n, k_1  , \ldots , k_r -1 )]  .\cr}
$$
>From this follows, upon summing, w.r.t. $k_1 , \ldots , k_r$, that
$$\displaylines{
a(n):= \sum_{ k_1 , \ldots , k_r } F( n, k_1 , \ldots , k_r ) ,\cr
\noalign{\hbox{satisfies that recurrence}}
S(N,n) a(n) \equiv 0  .\cr}
$$
However, using elimination is prohibitive. To find a {\eightrm
FAST} algorithm for multi-sum definite summation, we must first
find a multi-sum generalization of Gosper's algorithm. This
algorithm would input  a closed form $F( k_1 , \ldots , k_r )$ and
decide whether there exist closed form $G_i ( k_1 , \ldots , k_r
)$, ($i=1, \ldots , r$)
and find them in the affirmative case, such that:
$$
F  =  \sum_{i=1}^{r} \Delta_i G_i  .
$$
 
\medskip 
\centerline{\bf Epilogue written Feb. 1995 (for this version)}

\medskip 
The above was done in the following papers:
 
[WZ3] H.S. Wilf and D. Zeilberger,
{\it An algorithmic proof theory for hypergeometric
(ordinary and ``q") multisum/integral identities}, Invent. Math. 
{\bf 108}, 575-633 (1992).
 
[WZ4] H.S. Wilf and D. Zeilberger, {\it \hbox{\eightrm RATIONAL} function certification of
hypergeometric multi-integral/sum/"q" identities}, 
Bulletin of the Amer. Math. Soc. {\bf 27} 148-153 (1992).
 
[Z4] {\it Closed Form (pun intended!)}, 
in: "Special volume in memory of Emil Grosswald", M. Knopp and M. Sheingorn,
Contemporary Mathematics {\bf 143} 579-607, AMS, Providence (1993).
 
\medskip
Peter Paule found a great way to simplify the computer-generated
proofs for single-$q$-sums,
see:
 
[Pau] P. Paule, {\it Simple Computer Proofs for Rogers-Ramanujan type
Identities}, Elec. J. of Combinatorics {\bf 1}(1994), R10.
 
He and his students are currently developing farther ramifications as well
specialization and dualizations.
 
\medskip
Ira Gessel, has made a systematic study of 
specialization and dualization.
 
\medskip
My former student Sheldon Parnes has extended the algoritm for
`algebraic kernels', like in the generating function for the
Jacobi polynomials. See:
 
[EP] S.B. Ekhad and S. Parnes, {\it   A WZ-style proof of Jacobi polynomials' 
generating function. }, Discrete Mathematics 
{\bf 110}, 263-264 (1992).
 
[Par] S. Parnes, 
{\it A differential view of hypergeometric functions: algorithms and
implementation}, Ph.D. thesis,
Temple University, 1993.
[Available from University Microfilms, Ann Arbor, MI.]
 
\medskip
My student John Majewicz has extended Sister Celine's technique and
WZ-certification to Abel-type sums. See:
 
[EM] S.B. Ekhad and J.E. Majewicz, {\it A short WZ-style proof of Abel's
identity}, preprint, available by anon. ftp to {\tt ftp.math.temple.edu}
in file {\tt /pub/ekhad/abel.tex}.
 
[Ma] J.E. Majewicz, {\it WZ-style certification procedures and Sister's Celine's
technique for Abel-type sums},
preprint, available by anon. ftp to {\tt ftp.math.temple.edu}
in file {\tt /pub/jmaj/abel\_sum.tex}.
 
\medskip
Lily Yen, in a brilliant Penn thesis, under the direction of Herb
Wilf, has found effective a priori bounds for the number of
special cases one should check a given identity in order to
(rigorously!) know that it is true in general, in:
 
[Y] L. Yen, {\it Contributions to the proof theory of hypergeometric 
identities}, Ph.D. thesis, University of Pennsylvania, 1993.
[Available from University Microfilms, Ann Arbor, MI.]
 
\medskip
A beautiful exposition, as well as an {\eightrm AXIOM}
implementation, was written by Joachim Hornegger, in his
Erlangen {\it Diplomarbeit} under the direction of Volker Strehl:
 
[Ho] J. Hornegger, {\it Hypergeometrische Summation und polynomiale
Re\-kur\-sion}, Diplomarbeit, Erlangen, 1992.
 
\bye